\documentclass[12pt]{amsart}
\usepackage{amsmath,amssymb}
\newtheorem{theorem}{Theorem}
\newtheorem{proposition}[theorem]{Proposition}
\def\ds{\displaystyle}
\def\kob{\mbox{kob}}

\def\ve{\varepsilon}

\title{The Kobayashi balls of ($\mathbb{C}$-)convex domains}

\author{Nikolai Nikolov}

\author{Maria Trybu\l{}a}

\address{Institute of Mathematics and Informatics\\Bulgarian Academy
of Sciences\\ Acad. G. Bonchev 8, 1113 Sofia, Bulgaria\newline
\indent Faculty of Information Sciences\\
State University of Library Studies and Information Technologies\\
Shipchenski prohod 69A, 1574 Sofia, Bulgaria}\email{nik@math.bas.bg}
\address{Institute of Mathematics, Faculty of Mathematics and Computer
Science,\\ Jagiellonian University, \L ojasiewicza 6, 30-348
Krak\'ow, Poland} \email{maria.trybula@im.uj.edu.pl}



\thanks{The second named author is supported by
international PhD programme "Geometry and Topology in Physical
Models" of the Foundation for Polish Science. The paper was prepared
during her visit to the Institute of Mathematics and Informatics,
Bulgarian Academy of Science, October 2013 - April 2014.}

\begin{document}

\begin{abstract}  A pure geometric description of the Kobayashi balls
of ($\Bbb C$-)convex domains is given in terms of the so-called
minimal basis.
\end{abstract}

\maketitle

\section{Introduction and results}

Let $D$ be a domain in $\Bbb C^n$. Denote by $c_D$ and $l_D$ the
Carath\'eodory distance and the Lempert function of $D,$
respectively:
$$c_{D}(z,w)=\sup\{\tanh^{-1}|f(w)|:f\in\mathcal{O}(D,\mathbb{D}),\,
\hbox{ with }f(z)=0\},$$
$$l_{D}(z,w)=\inf\{\tanh^{-1}|\alpha|:\exists\varphi\in\mathcal{O}
(\mathbb{D},D) \hbox{ with }\varphi(0)=z,\varphi(\alpha)=w\},$$
where $\Bbb D$ is the unit disc. The Kobayashi distance $k_D$ is the
largest pseudodistance not exceeding $l_D.$

We are interested in a description of the Kobayashi balls near
boundary points of convex and, more generally, $\Bbb C$-convex
domains in terms of parameters that reflect the geometry of the
boundary. The first results in this direction can be found in
\cite[Theorems 1 and 5.1]{Aladro}, where the strongly pseudoconvex
case in $\Bbb C^n$ and the weakly pseudoconvex finite type case in
$\Bbb C^2$ are discussed with applications\footnote{See also
\cite{Sahin} for complex ellipsoids.} to invariant forms of
Fatou type theorems (for the boundary values). The weakly
pseudoconvex finite type case in $\Bbb C^2$, as well as the convex
finite type case in $\Bbb C^n,$ are treated in \cite[Propositions
8.8 and 8.9]{Mahajan} as byproducts of long considerations. The
strongly pseudoconvex case in $\Bbb C^n$ and the weakly pseudoconvex
finite type in $\Bbb C^2$ are particular cases of the pseudoconvex
Levi corank one case which is considered in \cite[Theorem
1.3]{Balakumar}. The behavior of the Kobayashi balls in all the
mentioned results is given in terms of the Levi geometry of the
boundary which is assumed smooth and bounded.

Our aim is to describe the Kobayashi balls of ($\Bbb C$-)convex
domains (not necessarily smooth and bounded) in terms of the
so-called minimal basis (cf. \cite{Conrad,Nikolov1,Nikolov4}. The
constants that appear depend only on the radius of the balls and the
dimension of the domains. The respective proof is short and pure
geometric. The obtained result covers \cite[Propositions 8.8 and
8.9]{Mahajan}.

Assume that $D$ contains no complex lines. Let $q\in D$ and
$d_D(q)=\mbox{dist}(q,\partial D).$ Choose $q^{1}\in\partial D$ so
that $\tau_1(q):=\lVert q^1-q\rVert=d_D(q).$ Put
$H_1=q+\textup{span}(q^{1}-q)^{\bot}$ and $D_1=D\cap H_1.$ Let
$q^2\in\partial D_1$ so that $\tau_2(q):=\lVert q^{2}-q\rVert=
d_{D_1}(q).$ Put $H_{2}=q+\textup{span}(q^{1}-q,q^{2}-q)^{\bot},$
$D_{2}=D\cap H_{2}$ and so on. Thus we get an orthonormal basis of
the vectors $\ds e_j=\frac{q^j-q}{\lVert q^{j}-q\rVert},$ $1\le j\le
n,$ which is called minimal for $D$ at $q,$ and positive numbers
$\tau_{1}(q)\le \tau_{2}(q)\le\dots\le \tau_{n}(q)$ (the basis and
the numbers are not uniquely determined). After rotation we may
assume that $e_1,e_2,\dots,e_n$ is the standard basis of $\Bbb C^n.$

Recall now that a open set $D$ in $\Bbb C^n$ is said to be (cf.
\cite{Anderson}):

$\bullet$ $\Bbb C$-convex if any non-empty intersection with a
complex line is a simply connected domain.

$\bullet$ linearly (weakly linearly convex) convex if for any
$a\in\Bbb C^n\setminus D$ ($p\in\partial D$) there exists a complex
hyperplane through $a$ which does not intersect $D.$

Note that convexity $\Rightarrow$ $\Bbb C$-convexity $\Rightarrow$
linear convexity $\Rightarrow$ weak linear convexity (cf.
\cite[Theorem 2.3.9 ii)]{Anderson} for the second implication).
Moreover, in the case of $C^1$-smooth bounded domains the last three
notions coincide (cf. \cite[Corollary 2.5.6]{Anderson}.

In view of this remark and the inequalities $c_D\le k_D\le l_D,$ we
have the following quantitative information about the
Carath\'eodory/Kobayashi/\ Lempert balls of ($\Bbb C$-)convex
domains.\footnote{By the Lempert theorem, $c_D=k_D=l_D$ in the
convex case, as well as in the bounded $C^2$-smooth $\Bbb
C$-convex case (cf. \cite{Nikolov3}).}

\begin{theorem}
Let $D$ be a domain in $\Bbb C^n,$ containing no complex lines, and
$q\in D.$ Assume that the standard basis of $\Bbb C^n$ is minimal
for $D$ at $q.$ Let $r>0.$

(i) If $D$ is weakly linearly convex, then
$$\begin{aligned}\max_{1\leq j\leq
n}\frac{|z_j-q_j|}{\tau_j(q)}<\frac{e^{2r}-1}{n(e^{2r}+1)}
&\Rightarrow\sum_{j=1}^n
\frac{|z_j-q_j|}{\tau_j(q)}<\frac{e^{2r}-1}{e^{2r}+1}\\
&\Rightarrow z\in D\mbox{ and }l_D(q,z)<r.\end{aligned}$$

(ii) If $D$ is convex, then $c_D(q,z)<r$ implies $\ds\max_{1\leq
j\leq n}\frac{|z_j-q_j|}{\tau_j(q)}< e^{2r}-1.$

(iii) If $D$ is $\Bbb C$-convex, then $c_D(q,z)<r$ implies
$\ds\max_{1\leq j\leq n}\frac{|z_j-q_j|}{\tau_j(q)}<e^{4r}-1.$
\end{theorem}

So there exist constants $c'=c'(r,n)$ and $c''=c''(r)$ such that
$$\Bbb D(q_1,c'\tau_1(q))\times\dots\times\Bbb
D(q_n,c'\tau_n(q))\subset\kob_D(q,r)$$
$$\subset\Bbb
D(q_1,c''\tau_1(q))\times\dots\times\Bbb D(q_n,c''\tau_n(q)),$$
where $\kob_D(q,r)$ is the Kobayashi ball $\{z\in D:k(q,z)<r\}$ and
$\Bbb D(p,r)=\{z\in\Bbb C:|z-p|<r\}.$ By \cite[Lemma 3.10]{Conrad},
the sizes of these polydiscs are comparable (in terms of small/big
constant depending on $D$) with the sizes of polydiscs in
\cite{Balakumar,Mahajan} arising from the Levi geometry of the
boundary. Thus Theorem 1 extends \cite[Propositions
8.9]{Mahajan}.

Note also that if $D$ is a proper $\Bbb C$-convex domain in $\Bbb
C^n$ containing complex line, then it is biholomorphic to
$D'\times\Bbb C^{n-k},$ where $D'$ is a bounded domain in $\Bbb
C^k,$ $0<k<n.$ (cf. Proposition 3 and the preceding remark in
\cite{Nikolov2}). So $\tau_k(q)<\infty=\tau_{k+1}(q)$ and it is easy
to see that Theorem 1 remains true.

To prove Theorem 1, we need the planar cases of the following

\begin{proposition}\label{2} (i) Let $D$ be proper convex
domain in $\Bbb C^n.$ Then \textup{(}cf.
\cite[(2)]{Nikolov5}\textup{)}
$$c_D(z,w)\ge\frac12\log\frac{d_D(z)}{d_D(w)}.$$
Moreover, if $n=1,$ then
$$c_D(z,w)\ge\frac12\log\left(1+\frac{|z-w|}{d_D(w)}\right).$$

(ii) Let $D$ be proper $\Bbb C$-convex domain in $\Bbb C^n.$ Then
$$c_D(z,w)\ge\frac14\log\frac{d_D(z)}{d_D(w)}.$$
Moreover, if $n=1,$ then
$$c_D(z,w)\ge\frac14\log\left(1+\frac{|z-w|}{d_D(w)}\right).$$
\end{proposition}

The constants $1/2$ and $1/4$ are sharp as the examples $D=\Bbb D$
and $D=\Bbb C_\ast\setminus\Bbb R^+$ show. Note that in the $\Bbb
C$-convex case the weaker estimate
$$c_D(z,w)\ge\frac14\log\frac{d_D(z)}{4d_D(w)}$$
is contained in \cite[Proposition 2]{Nikolov5}

Theorem 1 has a local version.

\begin{proposition}\label{local}
Let $D$ be a domain in $\Bbb C^n$ whose boundary contains no affine
discs through $a\in\partial D.$ Assume that the standard basis of
$\Bbb C^n$ is minimal for $D$ at $q\in D.$ Let $r>r'>0.$

(i) If $D$ is weakly linearly convex near $a$, then
$$\begin{aligned}\max_{1\leq j\leq
n}\frac{|z_j-q_j|}{\tau_j(q)}<\frac{e^{2r}-1}{n(e^{2r}+1)}
&\Rightarrow\sum_{j=1}^n
\frac{|z_j-q_j|}{\tau_j(q)}<\frac{e^{2r}-1}{e^{2r}+1}\\
&\Rightarrow z\in D\mbox{ and }l_D(q,z)<r.\end{aligned}$$ for $q$
sufficiently close to $a.$

(ii) If $D$ is convex near $a,$ then $k_D(q,z)<r'$ implies
$\ds\max_{1\leq j\leq n}\frac{|z_j-q_j|}{\tau_j(q)}<e^{2r}-1$ for
$q$ sufficiently close to $a.$

(iii) If $D$ is $\Bbb C$-convex near $a$ and bounded, then
$k_D(q,z)<r'$ implies $\ds\max_{1\leq j\leq
n}\frac{|z_j-q_j|}{\tau_j(q)}$ $<e^{4r}-1$ for $q$ sufficiently
close to $a.$
\end{proposition}

By any of the above three notions of convexity near $a$ we mean that
there exists a neighborhood $U$ of $a$ such that $D\cap U$ is an
open set with the respective global convexity.

Note that in the convex case, as well as in the $C^1$-smooth $\Bbb
C$-convex case, if $\partial D$ contains no affine discs through
$a,$ then $\partial D$ contains no analytic discs through $a$ (cf.
\cite[Propoisition 7]{Nikolov4}).

\section{Proofs}

\noindent{\it Proof of Theorem 1.} (i) Since $D$ contains
the discs $\Bbb D(q_1,\tau_1(q)),\dots,$ $\Bbb D(q_n,\tau_n(q))$
(lying in the respective coordinate complex planes), it contains
their convex hull
$$C=\{\zeta\in\Bbb C^n:h(\zeta)=\sum_{j=1}^n\frac{|\zeta_j-q_j|}{\tau_j(q)}<1\}$$
(cf. \cite[Lemma 15]{Nikolov4}). Then
$$l_D(q,z)\le l_C(q,z)=\tanh^{-1}h(z)$$
(cf. \cite[Proposition 3.1.10]{Jarnicki}) which implies (i).
\smallskip

Before proving (ii) and (iii) note that by ($\Bbb C$-)convexity and
the construction of the minimal basis there exists a complex
hyperplane $q^{j+1}+W_j$ through $q^{j+1}$ that is disjoint from
$D,$ $j=0,\dots,n-1.$ It is not difficult to see that $W_j$ is given
by the equation
$$\alpha_{j,1}\zeta_{1}+\cdots+\alpha_{j,j}\zeta_{j}+\zeta_{j+1}=0.$$
Let $\Lambda:\Bbb C^n\to\Bbb C^n$ be the linear mapping with matrix
whose rows are given by the vectors
$(\alpha_{j,1},\dots,\alpha_{j,j},1,0,$ $\dots,0).$ Set
$\Lambda_q(\zeta)=q+\Lambda(\zeta-q).$ Note that $G=\Lambda_q(D)$ is
a ($\Bbb C$-)convex domain. Denote by $G_j$ the projection of $G$
onto $j$-th coordinate plane. Then $G\subset G'=G_1\times\dots\times
G_n$ and the product formula for the Carath\'eodory distance (cf.
\cite[Theorem 9.5]{Jarnicki}) implies that
\begin{equation}\label{1}
c_D(q,z)\ge c_{G'}(q,\Lambda_q(z))=\max_{1\le j\le
n}c_{G_j}(q_j,z_j).
\end{equation}
Observe also that $d_{G_j}(q_j)=\tau_j(q).$
\smallskip

(ii) If $D$ is a convex domain, then $G_j$ is a convex domain.
Hence, by Proposition \ref{2} (i),
$$c_{G_j}(q_j,z_j)\ge\frac12\log\left(1+\frac{|z_j-q_j|}{\tau_j(q)}\right)$$
and (ii) follows from here and \eqref{1}.
\smallskip

(iii) If $D$ is a $\Bbb C$-convex domain, then $G_j$ is a simple
connected domain (cf. \cite[Theorem 2.3.6]{Anderson}). Hence, by
Proposition \ref{2} (ii),
$$c_{G_j}(q_j,z_j)\ge\frac14\log\left(1+\frac{|z_j-q_j|}{\tau_j(q)}\right)$$
and (iii) follows from here and \eqref{1}.
\smallskip

\noindent{\it Proof of Proposition \ref{2}.} After translation and
rotation, we may assume that $0\in\partial D$ and
$w=(d_D(w),0,\dots,0).$
\smallskip

(i) We have that $D\subset\Pi^+=\{\zeta\in\Bbb C^n:\mbox{Re
}\zeta_1>0\}$ and hence
$$c_D(z,w)\ge
c_{\Pi^+}(z,w)=\tanh^{-1}\left|\frac{z_1-w_1}{z_1+{\overline
w_1}}\right|$$
$$\ge\tanh^{-1}\frac{|z_1-w_1|}{|z_1-w_1|+2d_D(w)}=
\frac12\log\left(1+\frac{|z_1-w_1|}{d_D(w)}\right).$$

(ii) It follows by weak linear convexity that $D\cap\{\zeta_1\in\Bbb
C^n:\zeta_1=0\}=\varnothing.$ Denote by $D_1$ the projection of $D$
onto the $\zeta_1$-plane. Let $\gamma_G$ the Carath\'eodory metric
of a domain $G$ in $\Bbb C^k:$
$$\gamma_G(\zeta;X)=\sup\{|f'(\zeta)X|:f\in\mathcal O(G,\Bbb
D)\},\quad\zeta\in G,\ X\in\Bbb C^k.$$ The K\"{o}be $1/4$ theorem
implies that
$$\gamma_{D_1}(\zeta_1;e_1)\ge\frac{1}{4d_{D_1}(\zeta_1)}\ge\frac{1}{4|\zeta_1|}.$$
Since $D_1$ is a simply connected domain (cf. \cite[Theorem
2.3.6]{Anderson}), then
$$c_D(z,w)\ge
c_{D_1}(z_1,w_1)=\inf_s\int_0^1\gamma_{D_1}(s(t);s'(t)dt
\ge\frac14\inf_s\int_0^1\left|\frac{s'(t)}{s(t)}\right|dt,$$ where
the infimum is taken over all smooth curves $s:[0,1]\to D_1$ with
$s(0)=z_1$ and $s(1)=w_1$ (cf. \cite{Jarnicki}).

Set now
$$d(\zeta_1,\eta_1)=\log\max(1+|1-\zeta_1/\eta_1|,1+|1-\eta_1/\zeta_1|).$$
It is easy to check that $d$ is a distance on $\Bbb
C_\ast$\footnote{Let $a,b,c\in\Bbb C_\ast$ and $d_1=1-a/b,$
$d_2=1-b/c,$ $d_3=1-a/c.$ We may assume that $d(a,c)=\log(1+|d_3|).$
Then
$$d(a,b)+d(b,c)\ge\log(1+|d_1|)+\log(1+|d_2|)$$
$$=\log(1+|d_1|+|d_2|+|d_3-d_2-d_1|)\ge\log(1+|d_3|)=d(a,c).$$}
with ``derivative"
$$\lim_{\lambda\nrightarrow 0}
\frac{d(\zeta_1,\zeta_1+\lambda)}{|\lambda|}=\frac{1}{|\zeta_1|}.$$
Then (cf. \cite[Lemma 4.3.3) (d)]{Jarnicki})
$$\inf_s\int_0^1\left|\frac{s'(t)}{s(t)}\right|dt\ge d(z_1,w_1)$$
and hence
$$c_D(z,w)\ge\frac14d(z_1,w_1)\ge\frac14
\log\left(1+\frac{|z_1-w_1|}{d_D(w)}\right).$$

\noindent{\it Proof of Proposition \ref{local}.} (i) Using Theorem
1 (i), it is enough to show that
\begin{equation}\label{tau}
\lim_{q\to a}\tau_n(q)=0.
\end{equation}
Assume the contrary. Then there exists a sequence of
points $(q^j)\to a$ such that $(\tau_n(q^j))\to\ve>0$ and $(e^j)\to
e,$ where $e^j$ is the last vector of the minimal basis for $D$ at
$q^j.$ We may find a bounded neighborhood $U$ of $a$ such that
$D\cap U$ is a weakly linearly convex open set. Shrinking $\ve $ (if
necessary), it follows that the $e$-directional disc $\Delta$ with
center $q$ and radius $\ve$ is a limit of affine discs in $D\cap U.$
Since $D\cap U$ is a taut open set (cf. \cite[Proposition
1.5]{Nikolov3}), then $\Delta\subset\partial D,$ a contradiction.
\smallskip

(ii) Having in mind Theorem 1 (ii), it is enough to show
the following. \vskip5pt

{\it Claim 1.} Let $U$ be a neighborhood of $a$ such that $D\cap U$
is convex. There exist neighborhoods $W\subset V\subset U$ of $a$
such that if $q\in D\cap W$ and $k_D(q,z)<r',$ then $z\in V$ and
$pk_{D\cap U}(q,z)\le k_D(q,z),$ where $p=r'/r.$ \vskip5pt

To prove this claim, recall that $k_D$ is the integrated form of the
Kobayashi metric
$$\kappa_D(\zeta;X)=\inf\{|\alpha|:\exists\varphi\in\mathcal O(\Bbb D,D) \hbox{
with }\varphi(0)=z,\alpha\varphi'(0)=X\}$$ (cf. \cite{Jarnicki}).
Fix an $\ve>0.$ Then we may find a smooth curve $s:[0,1]\to D$ such
that $s(0)=q,$ $s(1)=z$ and
$$k_{D}(q,z)+\ve>\int_{0}^{1}\kappa_{D}(s(t);s'(t))dt.$$

Since $D\cap U$ is convex and its boundary contains no affine discs
through $a,$  then $a$ is a peak point for $D\cap U$ (cf.
\cite[Theorem 6]{Nikolov7}). Hence the strong localization property
for the Kobayashi metric holds (cf. \cite[Theorem 1 and Corollary
2]{Nikolov6}). So there exists a neighborhood $V\subset U$ of $a$
such that
$$\kappa_{D}(\zeta;X)\geq p\kappa_{D\cap U}(\zeta;X),
\quad\zeta\in D\cap V,\ X\in\Bbb C^n.$$

Set $t'=\sup\{t:s([0,t])\subset V\}$ and $z'=s(t').$ Then
\begin{multline*}
r'+\ve>k_{D}(q,z)+\ve>\int_0^{t'}\kappa_{D}(s(t);s'(t))dt\\
\geq p\int_0^{t'}\kappa_{D\cap U}(s(t);s'(t))dt \geq p k_{D\cap
U}(q,z')\geq pc_{D\cap U}(q,z').
\end{multline*}

Taking a peak function for $D\cap U$ at $a$ as a competitor in the
definition of $c_{D\cap U},$ it follows that
$$\lim_{q\to a}\inf_{\zeta\not\in V}c_{D\cap U}(q,\zeta)=+\infty.$$
Therefore, we may find a neighborhood $W\subset V$ such that if
$q\in W,$ then $z'\in V.$ Therefore $z'=z$ and the claim follows by
letting $\ve\to 0.$
\smallskip

(iii) In this case we do not know if $a$ is a local peak point.
However, Theorem 1 (iii) together with two small
modifications in the previous proof imply the
desired result.

1) The strong localization for the Kobayashi metric follows by the
property

$$\lim_{q\to a}\inf_{\zeta\not\in\mathcal U}l_D(q,\zeta)=+\infty,$$
where $\mathcal U$ is any neighborhood of $a$ (cf. \cite[Proposition
7.2.9]{Jarnicki}).

It is easy to see that this property is a consequence of the following
\smallskip

\textit{Claim 2.} If $(\varphi_j)\subset\mathcal{O}(\Bbb D,D)$ and
$\varphi_j(0)\to a,$ then $\varphi_j\rightrightarrows a.$
\smallskip

To prove Claim 2, assume the contrary. Since $D$ is bounded, then,
passing to a subsequence (if necessary), we may suppose that
$\varphi_j\rightrightarrows\varphi\in\mathcal
O(\mathbb{D},\overline{D})$ and $\varphi\neq a.$ Using again that
$D$ is bounded, we may find an $s\in(0,1)$ such that
$\varphi_j(s\Bbb D)\subset U$ for any $j.$ Note now that
\cite[Proposition 3]{Nikolov2} (see also \cite[Proposition
1.5]{Nikolov3}) implies the tautness of bounded $\mathbb{C}$-convex
domains. Then $D\cap U$ is a taut domain and hence $\varphi(s\Bbb
D)\in\partial D.$ Since $\partial D$ contains no affine discs
through $a\in\partial D,$ we get similarly to the proof of
\cite[Proposition 7]{Nikolov4} that $\varphi(s\Bbb D)=\{a\}.$ Then
the identity principle implies that $\varphi=a.$ This contradiction
completes the proof of Claim 2.
\smallskip

2) $\displaystyle \lim_{q\to a}\inf_{\zeta\not\in V}k_{D\cap
U}(q,\zeta)=+\infty.$\footnote{This equality instead of the same for 
$c_{D\cap U}$ can be used to finish the proof of Proposition \ref{local} (ii).}

This follows by Theorem 1 (iii) and the equality
\eqref{tau}.


\begin{thebibliography}{}

\bibitem{Anderson} M. Andersson, M. Passare, R. Sigurdsson,
{\it Complex convexity and analytic functionals}, Birkh\"auser,
Basel-Boston-Berlin, 2004.

\bibitem{Aladro} G. Aladro, {\it The comparability of the Kobayashi
approach region and the admissible approach region}, Illinois J.
Math. 33 (1989), 42-63.

\bibitem{Balakumar} G. P. Balakumar, P. Mahajan, K. Verma,
{\it Bounds for invariant distances on pseudoconvex Levi corank nne
domains and applications}, manuscripta math. (to appear);
arXiv:1303.3439.

\bibitem{Conrad} M. Conrad, {\it Nicht isotrope Absch\"atzungen f\"ur
lineal konvexe Gebiete endlichen Typs}, Dissertation, Universit\"at
Wuppertal, 2002.

\bibitem{Jarnicki} M. Jarnicki, P. Pflug, {\it Invariant distances
and metrics in complex analysis}, de Gruyter, Berlin-New York
(1993).

\bibitem{Mahajan} P. Mahajan, K. Verma, {\it Some aspects of the
Kobayashi and Carath\'eodory metrics on pseudoconvex domains}, J.
Geom. Anal. 22 (2012), 491-560.

\bibitem{Nikolov6} N. Nikolov, {\it Localization of invariant metrics},
Arch. Math. 79 (2002) 67-73.

\bibitem{Nikolov7} N. Nikolov, P. Pflug, {\it Behavior of the Bergman
kernel and metric near convex boundary points}, Proc. Amer. Math.
Soc. 131 (2003), 2097-2102.

\bibitem{Nikolov1} N. Nikolov, P. Pflug, {\it Estimates for the
Bergman kernel and metric of convex domains in $\Bbb C^n$}, Ann.
Polon. Math. 81 (2003), 73-78.

\bibitem{Nikolov2} N. Nikolov, P. Pflug, W. Zwonek,
{\it An example of a bounded $\Bbb C$-convex domain which is not
biholomorphic to a convex domain}, Math. Scand. 102 (2008), 149-155.

\bibitem{Nikolov3}  N. Nikolov, P. Pflug, P. J. Thomas, W. Zwonek,
{\it On a local characterization of pseudoconvex domains}, Indiana
Univ. Math. J. 58 (2009), 2661-2671.

\bibitem{Nikolov4} N. Nikolov, P. Pflug, W. Zwonek,
{\it Estimates for invariant metrics on  $\Bbb C^n$-convex domains},
Trans. Amer. Math. Soc. 363 (2011), 6245-6256.

\bibitem{Nikolov5} N. Nikolov, {\it Estimates of invariant metrics on
``convex domains''}, Ann. Mat. Pura Appl., DOI
10.1007/s10231-013-0345-7.

\bibitem{Sahin} S. Sahin, {\it Poletsky--Stessin Hardy spaces on complex
ellipsoids in $\Bbb C^n$}, arXiv:\\1403.2506.

\end{thebibliography}
\end{document}